**Linear-quadratic control of Volterra integral systems, and extensions.**

S. A. Belbas, Mathematics Department, University of Alabama, Tuscaloosa, AL 35487-0350, USA; e-mail sbelbas@gmail.com

W. H. Schmidt, Institut für Mathematik und Informatik, Universität Greifswald, Germany.

Abstract. We study linear-quadratic optimal control problems for Voterra systems, and Volterra optimal control problems that are linear-quadratic in the control but generally nonlinear in the state. In the case of linear-quadratic Volterra control, we obtain sharp necessary and sufficient conditions for optimality. For problems that are linear-quadratic in the control only, we obtain a novel form of necessary conditions in the form of double Volterra equations, for which we use the appellation "Riccati-Volterra equations"; we prove the solvability of such equations, which are the Volterra integral counterpart of the so-called "state-dependent Riccati equations".

1. Introduction.

This paper aims to present cases of optimal control problems for integral equations that are solvable with methods analogous to the Riccati equations of controlled differential systems. The well-known Riccati differential equation has the general form

$$\frac{dS(t)}{dt} = A(t)S(t) + S(t)A^T(t) + S(t)K(t)S(t) + \Phi(t) \tag{1.1}$$

and it arises the optimal control of a linear differential system with a cost functional that is quadratic in the state and the control. (Riccati equations arise also in other contexts.) The Riccati equation of linear-quadratic control of ODE systems is considered as one of the early success stories of optimal control and feedback systems. An up-to-date outline of Riccati equations, including an outline of its history, is contained in Bittanti et al. (1991) and Lewis (1992).

One of the issues we address here is the form of a counterpart of the Riccati equation for control systems governed by integral equations. It will turn out that an analogue of the state-dependent Riccati equation exists for Volterra control systems that are linear in the control, with a cost functional that is quadratic in the control, and it takes the form

$$z(t) = z_0(t) + \int_t^T [K(t,s)z(s) + z^T(s)L(s,t)]\, ds + \int_t^T \int_t^T z^T(s)M(s,t,\sigma)z(\sigma)\, ds\, d\sigma \tag{1.2}$$

where the arrays *K*, *L*, *M* also depend on the state of the system and on the final time *T*.
In the case of LQ controlled differential systems, an optimal control is obtained in the form of instantaneous linear feedback, $u^*(t) = G(t)x(t)$, where *G* is obtained from the system dynamics, the form of the cost functional, and the solution of the Riccati equation. In the case of Volterra control systems, an optimal control is also obtained in the form of anticipatory feedback (complementary concept to the concept of causal feedback),

$$u^*(t) = A_0(T,t)\, y(T) + A_1(T,t)(t)y(t) + \int_t^T A_2(T,t,s)y(s)\, ds \tag{1.3}$$

where $y(\cdot)$ is the state of the controlled Volterra system. In general, instantaneous feedback, and causal feedback, are not possible for the optimal control of Volterra systems. This observation applies to full Volterra equations of the form

$$y(t) = y_0(t) + \int_0^t [f_0(t,s,y(s)) + F_1(t,s,y(s))u(s)]\, ds\,. \tag{1.4}$$

It may be noted that integral formulations of differential systems, i.e. equations of the form

$$y(t) = y_0(t) + \int_0^t f(s,y(s),u(s))\, ds\,,$$

are also (sometimes) referred to as “Volterra integral equations”, and, for such memoryless “Volterra equations”, an optimal control is representable in the form of causal feedback (Pritchard and You (1996)).

Optimal control problems for Volterra integral equations have been treated in Bakke (1974), Schmidt (1980), Schmidt (1982), Belbas and Schmidt (2009). The references de Acutis (1985), Connor (1972), and Pandolfi (2018) deal with linear-quadratic problems for Volterra integral and integro-differential equations.

Our present work deals with more general problems for Volterra integral systems. We shall freely utilize, among other things, the necessary conditions for optimal control of Volterra equations Schmidt (1980) and Schmidt (1982).

## 2. Definitions and background information on second-degree functionals in infinite dimensional functional spaces.

In this section, we gather a few basic definitions and results pertaining to quadratic functionals on $L^2$ spaces.

We consider a functional of the form

$$E := \tfrac{1}{2}\iint_{G\times G} w^T(x)K_2(x,y)w(y)\,dx\,dy + \int_G \{\tfrac{1}{2}w^T(x)K_1(x)w(x) + r_0^T(x)w(x)\}\,dx \,. \tag{2.1}$$

Here, $G$ is a bounded open set in $\mathbb{R}^d$ , $w$ is an element of $L^2(G;n)$, the Hilbert space of square-integrable functions defined on $G$ with values in $\mathbb{R}^n$ .
The quadratic functional is well defined on $L^2(G;n)$ if the matrix-valued kernels $K_1$ and $K_2$ have square integrable entries and $r_0$ is in $L^2(G;n)$.
For additional purposes, over and above having well-defined quadratic forms, we may postulate stronger properties, such as continuity of $K_1, K_2$ on $\bar{G}, \bar{G}\times\bar{G}$, respectively.

We note that we may assume, without loss of generality, that the kernels $K_1$ and $K_2$ are <u>symmetric</u>, in the sense that $K_1^T(x) = K_1(x)$ and $K_2^T(x,y) = K_2(y,x)$.

We assume

<u>A.2.1</u>. The matrix-valued symmetric functions, with real entries, $K_1$ and $K_2$ are square-integrable functions in the sense that the Frobenius norms
$|K_1(x)|_F = \sqrt{\mathrm{tr}((K_1(x))^2)},\ |K_2(x,y)|_F = \sqrt{\mathrm{tr}((K_2(x,y))^2)}$ are square-integrable functions, and $r_0 \in L^2(G;n)$. (The prefix "tr" signifies the trace of a matrix.) ///

The absolute value symbol, $|\cdot|$, signifies the Euclidean vector norm of a finite-dimensional vector, or a matrix norm of a finite-dimensional matrix. We use the standard norm symbol, $\|\cdot\|$, for norms of vectors or operators in functional (infinite-dimensional) spaces.

Our definition of the Frobenius norms of matrices is adapted to symmetric matrices with real entries; a slightly different definition is standard for more general matrices.

We set

$$\| K_1 \|_F = \sqrt{\int_G \ | K_1(x) |_F^2 \ dx} \ , \| K_2 \|_F = \sqrt{\int_G \int_G \ | K_2(x,y) |_F^2 \ dx\, dy} \ ;$$

$$\| (K_1, K_2) \|_F := \| K_1 \|_F + \| K_2 \|_F \ .$$

We have the following inequality:

$$\left| \int_G \ v^T(x) K_1(x) w(x)\, dx + \int_G \int_G \ v^T(x) K_2(x,y) w(y)\, dx\, dy \right| \le \| (K_1, K_2) \|_F \| v \|_{L^2(G,n)} \| w \|_{L^2(G,n)}$$

The proof is, essentially, repeated application of the Cauchy-Schwarz inequality; we omit the details. ///

Definition 2.1. We shall say that the pair of matrix-valued kernels ($K_1$ , $K_2$) is a pair that generates a positive-definite integral form if $K_1(x)$ is invertible for all $x$ in $\overline{G}$ and, for every nonzero $w$ in $\left(L^2(G;n)\right)$ , we have

$$\tfrac{1}{2} \int_G \ w^T(x) K_1(x) w(x)\, dx + \tfrac{1}{2} \iint_{G \times G} \ w^T(x) K_2(x,y) w(y)\, dx\, dy > 0$$

When the strict inequality above is replaced by a non-strict inequality ( $\ge 0$ ), we shall say that the pair ($K_1$ , $K_2$) generates a nonnegative definite or positive semi-definite quadratic form. ///

Definition 2.2. We shall say that the matrix-valued function $K_1(x)$ is coercive on $G$ if there exists an $\alpha > 0$ such that, for almost all $x$ in $G$, we have $\xi^T K_1(x) \xi \ge \alpha |\xi|^2 \quad \forall \xi \in \mathbb{R}^n$. (The crucial condition here is the uniformity with respect to $x \in G$ , i.e. the coercivity constant $\alpha$ is independent of $x$.)
We shall say that $K_1(\cdot)$ is bounded on $G$ if there exists a constant $C$ such that, for almost all $x$ in $G$ , we have $|K_1(x)|_F \le C$ .
We shall say that the pair of kernels $(K_1, K_2)$ determines a coercive quadratic form, with coercivity constant $\beta > 0$ if, for every $w \in L^2(G;n)$,

$$\int_G \ w^T(x) K_1(x) w(x)\, dx + \iint_{G \times G} \ w^T(x) K_2(x,y) w(y)\, dx\, dy \ge \beta \| w \|_{L^2(G,n)}^2 .$$ ///

Remark 2.1. There is also a somewhat different definition of positive definiteness in contexts other than optimization and optimal control. We do not repeat that definition here, but we refer to [F, H, W]. Our definition of positive definiteness is consistent with the rubric "functions of positive type" in [M]. ///

Remark 2.2. In infinite dimensional spaces (like $L^2(G;n)$), coercive quadratic forms are a proper subset of positive definite quadratic forms. To justify this statement, it suffices to give examples of forms that are positive definite but not coercive.

Example for $K_1(\cdot)$: We take $K_1(x)=\alpha(x)I_n$, where $I_n$ is the $n\times n$ identity matrix, and $\alpha(\cdot)$ is real-valued, continuous and bounded on $G$, $\alpha(x)>0$ for all $x\in G\setminus\{x_0\}$, $\alpha(x_0)=0$ for a particular point $x_0\in G$. Then for every $\varepsilon>0$, we can find a $\delta>0$ such that $0<\alpha(x)<\varepsilon$ for all $x\in B(x_0,\delta)\setminus\{x_0\}$, where $B(x_0,\delta)$ is the open ball of radius $\delta$, centered at $x_0$, and $\delta>0$ is so small that $\overline{B(x_0,\delta)}\subseteq G$. We take a nonzero $w\in L^2(G;n)$ supported on $\overline{B(x_0,\delta)}$. Then $0<\int w^T(x)K_1(x)w(x)\,dx\le\varepsilon\int_G |w(x)|^2\,dx$. It is plain that the quadratic form defined by this kernel is positive definite but not coercive.

Example for $K_2(\cdot,\cdot)$: We take a complete orthonormal basis for $L^2(G;n)$, say $\{e_k : k\in\mathbb{N}\}$. We define a kernel

$$K_2(x,y)=\sum_{k=1}^{\infty}\lambda_k\,e_k(x)e_k^T(y)$$

where the series is convergent in the Frobenius norm and $\{\lambda_k : k\in\mathbb{N}\}$ is a sequence of positive terms with $\lim_{k\to\infty}\lambda_k=0$. We take $w_k=e_k$. Then

$Q_2(w_k)\equiv\int_G\int_G w^T(x)K_2(x,y)w(y)\,dx\,dy=\lambda_k=\lambda_k\|w_k\|^2_{L^2(G;n)}$, thus $Q_2$ is not coercive. On the other hand, for a general $w\in L^2(G;n)$, we have

$$Q_2(w)=\sum_{k=1}^{\infty}\lambda_k\left(\int_G w^T(x)e_k(x)\,dx\right)^2$$

which shows the positive definiteness of $Q_2$. ///

We invoke the results of (Belbas and Schmidt 2016) to state:

Theorem 2.1. When $K_1(\cdot)$ is bounded and coercive, and the pair ($K_1$, $K_2$) is a pair that generates a positive-definite integral form, then:

(i) The Fredholm integral equation

$$K_1(x)w(x)+\int_G K_2(x,y)w(y)\,dy+r_0(x)=0$$

has a unique solution $w^*$;

(ii) The unique solution $w^*$ minimizes $E(w)$ over all $w$ in $L^2(G;n)$ . ///

3. <u>Linear Volterra control systems with general quadratic cost functional.</u>

We consider a linear controlled Volterra equation

$$y(t) = y_0(t) + \int_0^t [A(t,s)y(s) + B(t,s)u(s)] \tag{3.1}$$

with general quadratic cost functional

$$J = \tfrac{1}{2}y^T(T)P_0y(T) + \int_0^T [\tfrac{1}{2}y^T(t)P_1(t)y(t) + y^T(t)Q_1(t)u(t) + \tfrac{1}{2}u^T(t)R_1(t)u(t)]\,dt +$$
$$+\int_0^T q_0^T(t)u(t)\,dt + \int_0^T\int_0^T [\tfrac{1}{2}y^T(t)P_2(t,\tau)y(\tau) + y^T(t)Q_2(t,\tau)u(\tau) + \tfrac{1}{2}u^T(t)R_2(t,\tau)u(\tau)]\,d\tau\,dt \tag{3.2}$$

The matrix valued functions $A, B, P_i, Q_i, R_i, i=1,2$, and the vector-valued function $q_0$, are continuous in the relevant closed domains, i.e. $A$ and $B$ are continuous on $\bar{D}_1 = \{(t,s): 0 \le s \le t \le T\}$, $y_0, P_1, Q_1, R_1$ are continuous on $[0,T]$, $P_2, Q_2, R_2$ are continuous on $[0,T]\times[0,T]$.

The solution of the equation of state dynamics is representable as

$$y(t) = y_1(t) + \int_0^t B_1(t,s)u(s)\,ds \tag{3.3}$$

Indeed, if $S(t,s)$ is a resolvent kernel corresponding to direct kernel $A$, then

$$y_1(t) = y_0(t) + \int_0^t S(t,\sigma)y_0(\sigma)\,d\sigma\ ,\quad B_1(t,s) = B(t,s) + \int_s^t S(t,\sigma)B(\sigma,s)\,d\sigma \tag{3.4}$$

By using the representation (3.3) of the state, the cost functional is expressed in the form

$$J = \Phi(t) + \tfrac{1}{2}\int_0^T u^T(t)R_1(t)u(t)\,dt + \tfrac{1}{2}\int_0^T\int_0^T u^T(t)K_2(t,\tau)u(\tau)\,d\tau\,dt + \int_0^T q^T(t)u(t)\,dt \tag{3.5}$$

where $\Phi(t)$ is independent of the control, and

$$K_2(t,\tau) = B_1^T(T,t)P_0B_1(T,\tau) + \int_{t\vee\tau} B_1^T(s,t)P_1(s)B_1(s,\tau)\,ds +$$
$$+B_1^T(\tau,t)Q_1(\tau) + Q_1^T(t)B_1(t,\tau) + \int_t^T\int_\tau^T B_1^T(s,t)P_2(s,\sigma)B_1(\sigma,\tau)\,d\sigma\,ds + \tag{3.6}$$
$$+\int_t^T B_1^T(s,t)Q_2(s,\tau)\,ds + \int_\tau^T Q_2^T(\sigma,t)B_1(\sigma,\tau)\,d\sigma + R_2(t,\tau)$$

$$q^T(t) = y_1^T(T)P_0B_1(T,t) + \int_t^T y_1^T(s)P_1(s)B_1(s,t)\,ds + y_1^T(t)Q_1(t) + \\ + \int_0^T\int_t^T y_1^T(\tau)P_2(\tau,s)B_1(s,t)\,ds\,d\tau + \int_0^T y_1^T(\tau)Q_2(\tau,t)\,d\tau + q_0^T(t) \tag{3.7}$$

An optimal control $u^*$ must satisfy

$$R_1(t)u^*(t) + \int_0^T K_2(t,\tau)u^*(\tau)\,d\tau + q(t) = 0 \tag{3.8}$$

We assume:

A.0. $Q_1(t) = Q_2(t,\tau) = 0 \quad \forall(t,\tau)\in[0,T]\times[0.T]$.
A.1. $y_0, q_0, P_1, R_1$ are continuous on $[0,T]$; $A,B$ are continuous on $D=\{(s,t):0\le s\le t\le T\}$; $P_2$ , $R_2$ are continuous on $[0,T]\times[0,T]$.
A.2. $P_0$ is nonnegative definite; $R_1$ is positive definite on $[0,T]$; the pair $(P_1,P_2)$ generates a nonnegative quadratic form on $L^2(0,T;n)$; the pair $(R_1,R_2)$ generates a positive-definite quadratic form on $L^2(0,T;n)$ .

Theorem 3.1. Under assumptions A.0 through A.2, the integral equation (3.8) has a unique solution which is the wanted optimal control.

Proof: Condition (A.1) implies that $B_1(t,s)$ is well-defined and continuous on *D*, by standard results on linear Volterra integral equations.
It suffices to establish that the pair $(R_1,K_2)$ generates a nonnegative definite quadratic form on $L^2(0,T;n)$. Concerning the terms, on the RHS of (3.6) that contain $P_1$ and $P_2$, we observe that, for every $w\in L^2(0,T;n)$, we have

$$\int_0^T\int_0^T\int_{t\vee\tau}^T w^T(t)B_1^T(s,t)P_1(s)B_1(s,\tau)w(\tau)\,ds\,dt\,d\tau + \\ + \int_0^T\int_0^T\int_t^T\int_\tau^t w^T(t)B_1^T(s,t)P_2(s,\sigma)B_1(\sigma,\tau)w(\tau)\,ds\,d\sigma\,dt\,d\tau = \\ = \int_0^T \left(\int_0^s B_1(s,t)w(t)\,dt\right)^T P_1(s)\left(\int_0^s B_1(s,\tau)w(\tau)\,d\tau\right)ds + \\ + \int_0^T\int_0^T \left(\int_0^s B_1(s,t)w(t)\,dt\right)^T P_2(s,\sigma)\left(\int_0^\sigma B_1(\sigma,\tau)w(\tau)\,d\tau\right)ds\,d\sigma \ge 0 \tag{3.9}$$

An analogous estimate holds for the terms containing $P_0$.

The remaining conditions in (A.2) then imply that the pair of kernels $(R_1, K_2)$ fulfill the requirements of section 2 for existence and uniqueness of solution of (3.9) and the minimization property of that solution. ///

Remark 3.1. The assumption A.0 is not essentially restrictive, it is a condition ordinarily included in treatments of LQ control problems. The removal of the assumption A.0 would require amending the remaining assumptions with A.3 below:

A.3. The quadratic form determined by $(R_1, R_2)$ is coercive with sufficiently large coercivity constant.

The necessary magnitude of this coercivity constant depends on the Frobenius norms of the kernels containing $Q_1, Q_2$ on the RHS of (3.6). We omit the detailed calculation of this necessary magnitude of the coercivity constant, as it is straightforward but tedious, and not very informative in the present context. ///

4. Single-integral cost function and (non-causal) optimal feedback.

We consider the linear control system (3.), with cost functional

$$J = \tfrac{1}{2} y^T(T) P_0 y(T) + \int_0^T [\tfrac{1}{2} y^T(t) P_1(t) y(t) + y^T(t) Q(t) u(t) + \tfrac{1}{2} u^T(t) R(t) u(t)]\, dt \tag{4.1}$$

By analogy with the case of controlled ordinary differential equations (ODEs), we would anticipate an optimal control in the form of linear feedback. In the case of controlled ODEs, an optimal control is obtained as instantaneous linear feedback, i.e. the current value of the control is a linear function of the current value of the state. For our problem of controlled Volterra integral equations, it turns out that an optimal control is an anti-causal function of the state, of the form

$$u^*(t) = A_0(T,t) y(T) + A_1(T,t) y(t) + \int_t^T A_2(T,t,s) y(s)\, ds \tag{4.2}$$

This is the best possible result for LQC of Volterra control systems; instantaneous feedback, or causal feedback, are generally impossible.
We proceed to justify our claims. We shall utilize the Hamiltonian formalism, with co-state $\psi(t)$ , whose values are co-vectors (row vectors) of the same dimension as the state. The Hamiltonian is

$$\begin{aligned} H(t,T,y,Y,u,\psi(\cdot)) &= Y^T P_0 [A(T,t) y + B(T,t) u] + \tfrac{1}{2} y^T P_1(t) y + y^T Q_1(t) u + \tfrac{1}{2} u^T R(t) u + \\ &+ \int_t^T \psi(s) [A(s,t) y + B(s,t) u]\, ds \end{aligned} \tag{4.3}$$

The variable $Y$ occupies the slot of $y(T)$.

An optimal control satisfies

$$\nabla_u H = 0 \tag{4.4}$$

thus an optimal control $u^*$ is given by

$$u^{*T}(t) = -\left[ Y^T P_0 B(T,t) + y^T(t) Q_1(t) + \int_t^T \psi(s) B(s,t)\, ds \right] R_1^{-1}(t) \tag{4.5}$$

The co-state solves the Hamiltonian equation

$$\psi(t) = \nabla_y H(t, T, ...) \tag{4.6}$$

at is

$$\psi(t) = Y^T P_0 A(T,t) + y^T(t) P_1(t) + u^T(t) Q_1^T(t) + \int_t^T \psi(s) A(s,t)\, ds \tag{4.7}$$

By substituting the $u^*$ (as expressed by (4.6)) into (4.7), we obtain the linear Volterra equation for the co-state that corresponds to $u^*$ :

$$\begin{aligned} \psi^*(t) &= Y^T P_0 (A(T,t) - B(T,t) R_1^{-1}(t) Q_1^T(t)) + y^T(t)(P_1(t) - Q_1(t) R_1^{-1}(t) Q_1^T(t)) + \\ &+ \int_t^T \psi^*(s)(A(s,t) - B(s,t) R_1^{-1}(t) Q_1^T(t))\, ds \end{aligned} \tag{4.8}$$

We rewrite (4.8) in succinct form:

$$\psi^*(t) = Y^T M_0(T,t) + y^T(t) M_1(t) + \int_t^T \psi^*(s) N(s,t)\, ds \tag{4.9}$$

This is a Volterra integral equation in backward time. It has a resolvent kernel $\Sigma(\sigma,t)$ so that the solution is

$$\begin{aligned} \psi^*(t) &= Y^T M_0(T,t) + y^T(t) M_1(t) + \\ &+ \int_t^T [Y^T M_0(T,\sigma) + y^T(\sigma) M_1(\sigma)] \Sigma(\sigma,t)\, d\sigma \end{aligned} \tag{4.10}$$

Consequently (by substituting (4.10) into (4.5)),

$$\begin{aligned} u^{*T}(t) &= -[Y^T P_0 B(T,t) + y^T(t) Q_1(t)] R_1^{-1}(t) + \int_t^T [Y^T M_0(T,s) + y^T(s) M_1(s)] R_1^{-1}(t)\, ds + \\ &+ \int_t^T \int_t^s [Y^T M_0(T,s) + y^T(s) M_1(s)] \Sigma(s,\sigma) B(\sigma,t)\, d\sigma\, ds \end{aligned} \tag{4.11}$$

which, after transposition, gives $u^*(t)$ in the wanted feedback form.

In turn, substitution of the feedback form of $u^*$ into (4.1) gives the following Fredholm equation for the optimal trajectory:

$$\begin{aligned} y^*(t) &= y_1(t) + \left( \int_0^t B_1(t,s) A_0(T,s)\, ds \right) y^*(T) + \int_0^t B_1(t,s) A_1(T,s) y^*(s)\, ds + \\ &+ \int_0^T \int_0^{s \wedge t} B_1(t,\sigma) A_2(T,\sigma,s) y^*(s)\, d\sigma\, ds \end{aligned} \tag{4.12}$$

where $s \wedge t = \min(s,t)$.
In deriving the above equation, we have used a few simple calculus manipulations, which we do not show in detail.

Also, for $t = T$, we obtain the equation

$$
\begin{aligned}
y^*(T) &= y_1(T) + \left( \int_0^T B_1(T,s) A_0(T,s)\, ds \right) y^*(T) + \\
&+ \int_0^T \left[ B_1(T,s) A_1(T,s) + \int_0^s B_1(T,\sigma) A_2(T,\sigma,s)\, d\sigma \right] y^*(s)\, ds
\end{aligned}
\tag{4.13}
$$

Assuming invertibility of $I_n - \int_0^T B_1(T,s) A_0(T,s)\, ds$, we have

$$
\begin{aligned}
y^*(T) &= \left[ I_n - \int_0^T B_1(T,s) A_0(T,s)\, ds \right]^{-1} \cdot \\
&\cdot \left[ y_1(T) + \int_0^T \left[ B_1(T,s) A_1(T,s) + \int_0^s B_1(t,\sigma) A_2(T,\sigma,s)\, ds \right] y^*(s)\, ds \right]
\end{aligned}
\tag{4.14}
$$

Substitution of the above expression for $y^*(T)$ into (4.12) yields a Fredholm integral equation for the optimal trajectory. We omit the explicit form of this Fredholm integral equation.

## 5. State-dependent Riccati-Volterra integral equation.

State-dependent Riccati differential equations (for systems governed by ordinary differential equations) are well-known [BLT]. In this section we establish an analogous type of equations for Volterra control systems. This mathematically interesting, a it shows clearly the conceptual relationship between Riccati differential systems (the state-dependent variant thereof) to Riccati-Volterra equations.

We consider the controlled Volterra integral system that is affine in the control but generally nonlinear in the state:

$$y(t) = y_0(t) + \int_0^t [f_0(t,s,y(s)) + F_1(t,s,y(s))u(s)]\,ds \quad \left(0 < t \le T\right) \tag{5.1}$$

with a cost functional that is quadratic in the control:

$$J = \varphi_0(T, y(T)) + \int_0^T [g_0(t, y(t)) + g_1^T(t, y(t))u(t) + \tfrac{1}{2}u^T(t)G_2(t, y(t))u(t)]\,dt\,. \tag{5.2}$$

The Hamiltonian is

$$\begin{aligned} H(t,T,y,Y,u,\psi(\cdot)) &= \left(\nabla_Y\varphi_0(T,Y)\right)[f_0(T,t,y) + F_1(T,t,y)u] + \\ &+ g_0(t,y) + g_1^T(t,y)u + \tfrac{1}{2}u^T G_2(t,y)u + \int_t^T \psi(s)[f_0(s,t,y) + F_1(s,t,y)u]\,ds \end{aligned} \tag{5.3}$$

We assume positive definiteness of $G_2(t,y)$.

An optimal control satisfies

$$\nabla_u H = 0 \tag{5.4}$$

which leads to

$$u^{*T}(t) = -\left[ \left(\nabla_Y\varphi_0\right) F_1(T,t,y) + g_1^T(t,y) + \int_t^T \psi(s)F_1(s,t,y)\,ds \right] G_2^{-1}(t,y) \tag{5.5}$$

The Hamiltonian equation for the co-state is

$$\psi(t) = \nabla_y H \tag{5.6}$$

and in extensive form

$$\psi(t) = (\nabla_Y \varphi_0)\big(\nabla_y f_0(T,t,y)\big) - (\nabla_Y \varphi_0)\big(\nabla_y F_1(T,t,y)\big)u + \nabla_y g_0 + \big(\nabla_y g_1^T\big)u + \\ + \tfrac{1}{2}u^T\big(\nabla_y G_2\big)u + \int_t^T \psi(s)\big[\nabla_y f_0(s,t,y) + \nabla_y F_1(s,t,y)u\big]ds \tag{5.7}$$

(An explanatory notice may be appropriate here. We follow standard notational conventions. We use subscripts for the components of a vector, and superscripts for the components of a co-vector. The gradient of a real-valued function $f$ is a co-vector, $\big(\nabla_y f\big)^i = \dfrac{\partial f}{\partial y_i}$. For a vector-valued function $f = \big(f_j\big)$, $\big(\nabla_y f\big)^i_j = \dfrac{\partial f_j}{\partial y_i}$. For a matrix-valued function $M(y) = \big(M_j^k(y)\big)$, $\big(\nabla_y M\big)^{ki}_j = \left(\dfrac{\partial M_j^k}{\partial y_i}\right)$. For a quadratic expression, like $u^T\big(\nabla_y G_2\big)u$, the notation has the meaning $\big(u^T\big(\nabla_y G_2\big)u\big)^i = \sum_{j,k} u^j \dfrac{\partial\big(G_2\big)^k_j}{\partial y_i} u_k$. In this notation, $u = [u_j : 1 \le j \le m]$, $u^T = [u^k : 1 \le k \le m]$, and $u^k = \sum_j \delta^{jk} u_j (= u_k)$, $u_j = \sum_k \delta_{jk} u^k (= u^j)$, where the deltas are Kronecker's deltas. Our notation is a particular case of, and consistent with, the general notation of tensor analysis and Riemannian geometry.)

The co-state $\psi^*$ that corresponds to an optimal control $u^*$ satisfies an equation that results from substituting (5.7) into (5.5):

$$\psi^*(t) = \left(\nabla_Y \varphi_0\right)\Big[ \ \left(\nabla_y f_0(T,t,y)\right) - \left(\nabla_y F_1(T,t,y)\right) G_2^{-1}[F_1^T(t,t,y)\left(\nabla_Y \varphi_0\right)^T + g_1 \ \Big] +$$
$$+(-1)\left(\nabla_Y \varphi_0\right)\left(\nabla_y F_1(T,t,y)\right) G_2^{-1} \int_t^T \ F_1^T(s,t,y)\psi^{*T}(s)\, ds + \nabla_y g_0 - \left(\nabla_y g_1^T\right) G_2^{-1}[F_1^T(T,\ldots)\left(\nabla_Y \varphi_0\right)^T + g_1] +$$
$$+(-1)\left(\nabla_y g_1^T\right) G_2^{-1} \int_t^T \ F_1^T(s,t,y)\psi^{*T}(s)\, ds +$$
$$+\tfrac{1}{2}[\left(\nabla_Y \varphi_0\right) F_1(T,t,y) + g_1^T] G_2^{-1} \left(\nabla_y G_2\right) G_2^{-1}[F_1^T(T,\ldots)\left(\nabla_Y \varphi_0\right)^T + g_1] +$$
$$+\tfrac{1}{2}\int_t^T \ \psi^*(s) F_1(s,t,y)\, ds\, G_2^{-1} \left(\nabla_y G_2\right) G_2^{-1}[F_1^T(T,\ldots)\left(\nabla_Y \varphi_0\right)^T + g_1] +$$
$$+\tfrac{1}{2}[\left(\nabla_Y \varphi_0\right) F_1(T,t,y) + g_1^T] G_2^{-1} \left(\nabla_y G_2\right) G_2^{-1} \int_t^T \ F_1^T(s,t,y)\psi^{*T}(s)\, ds +$$
$$+\tfrac{1}{2}\int_t^T \int_t^T \ \psi^*(s) F_1(s,t,y) G_2^{-1} \left(\nabla_y G_2\right) G_2^{-1} F_1^T(\sigma,t,y)\psi^{*T}(\sigma)\, d\sigma\, ds \qquad (5.8)$$

This is the promised state-dependent Riccati-Volterra integral equation.

6. Solvability of Riccati-Volterra integral equations.

In this section, we prove sufficient conditions for existence and uniqueness of solution of a Riccati-Volterra equation of the form

$$\psi(t) = \psi_0(t) + \int_t^T \psi(s)K(s,t)\,ds + \tfrac{1}{2}\int_t^T \int_t^T \psi(s)\,L(s,t,\sigma)\,\psi^T(\sigma)\,ds\,d\sigma \tag{6.1}$$

The term $L(s,t,\sigma)$ is a function that takes values in the space of 3-index arrays,

$$L(s,t,\sigma) = \left[ L^i_{jk}(s,t,\sigma) : 1 \le i \le n, 1 \le j \le n, 1 \le k \le n \right]$$

and

$$\left[ \psi(s)\,L(s,t,\sigma)\,\psi^T(\sigma) \right]^i = \sum_{j,k} \psi^j(s) L^i_{jk}(s,t,\sigma) \psi^k(\sigma)$$

It is convenient to utilize a change of variables:

$$\begin{aligned} &\tilde t := T - t, \tilde s := T - s, \tilde\sigma := T - \sigma, \tilde\psi(\tilde t) := \psi(T - \tilde t), \tilde K(\tilde s, \tilde t) := K(T - \tilde s, T - \tilde t), \\ &\tilde L(\tilde s, \tilde t, \tilde\sigma) := L(T - \tilde s, T - \tilde t, T - \tilde\sigma), \tilde\psi_0(\tilde t) := \psi_0(T - \tilde t) \end{aligned} \tag{6.2}$$

Then the integral equation (6.1) becomes an equation in forward time:

$$\tilde\psi(\tilde t) = \tilde\psi_0(\tilde t) + \int_0^{\tilde t} \tilde\psi(\tilde s)\,\tilde K(\tilde s, \tilde t)\,d\tilde s + \tfrac{1}{2}\int_0^{\tilde t}\int_0^{\tilde t} \tilde\psi(\tilde s)\,\tilde L(\tilde s, \tilde t, \tilde\sigma)\,\tilde\psi^T(\tilde\sigma)\,d\tilde s\,d\tilde\sigma \tag{6.3}$$

We will present the results for integral equations in forward time, and for the scalar case. Thus we consider an equation

$$z(t) = z_0(t) + \int_0^t a(t,s)\,z(s)\,ds + \tfrac{1}{2}\int_0^t\int_0^t b(t,s,\sigma)\,z(s)\,z(\sigma)\,ds\,d\sigma \tag{6.4}$$

The real-valued functions $a$ and $b$ are continuous on $D_1 := \{(t,s) : 0 \le s \le t \le T\}$ and $D_2 := \{(t,s,\sigma) : 0 \le \min(s,\sigma),\ \max(s,\sigma) \le t \le T\}$ , respectively. We set

$$\bar a := \max_{(t,s)\in D_1} |\,a(t,s)\,|, \quad \bar b := \max_{(t,s,\sigma)\in D_2} |\,b(t,s,\sigma)\,|. \tag{6.5}$$

Every proof will require an *a priori* bound on the solution, and consequently some restrictions on the magnitude of *a* and *b*, and the time-horizon *T*.

We present one of the possible proofs. We seek solutions in the space $C(0,T;\mathbb{R})$ of continuous functions from $[0,T]$ into $\mathbb{R}$. We shall utilize the norm

$$\| z \|_{\mu} := \max_{0\le t\le T} e^{-\mu t} \, | z(t) |$$

and also

$$\| z \|_{0,s;\mu} := \max_{0\le \sigma\le s} e^{-\mu\sigma} \, | z(\sigma) | .$$

We denote by **S** the operator that represents the right-hand side of the integral equation:

$$(\mathbf{S}z)(t) := z_0(t) + \int_0^t a(t,s)\, z(s)\, ds + \tfrac{1}{2}\int_0^t\int_0^t b(t,s,\sigma)\, z(s)\, z(\sigma)\, ds\, d\sigma \tag{6.6}$$

We have

$$e^{-\mu t}\, | (\mathbf{S}z)(t) - z_0(t) | \le \int_0^t \overline{a}\, e^{\mu(s-t)} \| z \|_{0,s;\mu}\, ds + \int_0^t\int_0^t \overline{b}\, e^{\mu(s+\sigma-t)} \| z \|_{0,s;\mu} \| z \|_{0,\sigma;\mu}\, ds\, d\sigma \le$$
$$\le \overline{a}\, \frac{1-e^{-\mu t}}{\mu} \| z \|_{0,t;\mu} + \overline{b}\, \frac{\cosh t - 1}{\mu^2} \| z \|_{0,t;\mu}^2$$

Consequently

$$\| \mathbf{S}z \|_{\mu} \le \| z_0 \|_{\mu} + \overline{a}\, \frac{1-e^{-\mu T}}{\mu} \| z \|_{\mu} + \overline{b}\, \frac{\cosh(\mu T) - 1}{\mu^2} \| z \|_{\mu}^2 \tag{6.7}$$

We want a constant $C$ such that the condition $\| z \|_{\mu} \le C$ implies $\| \mathbf{S}z \|_{\mu} \le C$ .

We introduce the following notation:

$$c_0 = \| z_0 \|_{\mu} \,, \quad \alpha = \alpha(\mu,T) = \overline{a}\, \frac{1-e^{-\mu T}}{\mu} \,, \quad \beta = \beta(\mu,T) = \overline{b}\, \frac{\cosh(\mu T) - 1}{\mu^2} \tag{6.8}$$

We define the trinomial $\varphi$ by

$$\varphi(\xi) = \beta\xi^2 + (\alpha - 1)\,\xi + c_0 \,. \tag{6.9}$$

In this formulation, the problem has 5 parameters: $c_0, T, \overline{a}, \overline{b}, \mu$. (The other parameters that will appear in the proof are expressed in terms of these 5.)

We shall prove that, under certain conditions on these 5 parameters, the integral equation has a unique solution.

Lemma 6.1. We assume that the parameters $c_0, T, \overline{a}, \overline{b}, \mu$ satisfy

$\alpha < 1$, $c_0 \beta \le \frac{1}{4}(\alpha - 1)^2$, $c_0 \le \overline{r}$ (where $\overline{r}$ is the biggest real root of $\varphi$ ).

Then the set

$M := \{z \in C(0,T;\mathbb{R}) : \| z \|_\mu \le \overline{r}\}$
remains invariant under the operator **S**.

Proof: We are looking for a constant $C$ such that the inequality

$\| z \|_\mu \le C$

implies

$\| \mathbf{S}z \|_\mu \le C$ .

The trinomial $\varphi(\xi)$ satisfies $\varphi(\xi) \le 0$ for $\underline{r} \le \xi \le \overline{r}$ , where $\underline{r}, \overline{r}$ are the roots of $\varphi(\xi)$ (it is plain that, by dint of our stated assumptions, $\varphi(\xi)$ has real positive roots). Thus, for $\underline{r} \le \| z \| \le \overline{r}$ , we have

$\| \mathbf{S}z \|_\mu \le \| z \|_\mu$

thus (since $\| z \|_\mu \le \overline{r}$ )

$\| \mathbf{S}z \|_\mu \le \overline{r}$ .

On the other hand, for $0 \le \xi \le \underline{r}$, the trinomial $\varphi(\xi)$ s a decreasing function of $\xi$ , thus

$\varphi(\xi) \le \varphi(0) = \| z_0 \|_\mu \le \underline{r}$ .

Thus the implication

$0 \le \| z \|_\mu \le \overline{r} \;\; \Rightarrow \;\; \| \mathbf{S}z \|_\mu \le \overline{r}$

is valid. ///

Remark 3.1. The set of quintuples of parameters $(c_0, \mu, \bar{a}, \bar{b}, T)$ that satisfy the assumptions of the above lemma, is nonempty. This can be verified by noticing that $\alpha(\mu) \to 0$ as $\mu \to \infty$, so that $\alpha(\mu) < 1$ for $\mu$ sufficiently large, and the other assumptions are satisfied for sufficiently small values of the remaining parameters. ///

We have:

Theorem 6.1. Under the conditions of the lemma above, the Riccati-Volterra equation (6.4) has a unique solution.

Proof: Let $C$ be a constant such that $\| z \|_\mu \le C \;\Rightarrow\; \| \mathbf{S} z \|_\mu \le C$. (Such a constant $C$ exists, by lemma 6.1.) We set

$$M = \{ z \in C(0,T;\mathbb{R}) : \| z \|_\mu \le C \}.$$

Then $\mathbf{S}$ maps $M$ into itself. Further, for every two $z_1$, $z_2$ in $M$, we have

$$| b(t,s,\sigma)(z_1(s) z_1(\sigma) - z_2(s) z_2(\sigma) | \le \bar{b} [ | z_1(\sigma) | \, | z_1(s) - z_2(s) | + | z_2(s) | \, | z_1(\sigma) - z_2(\sigma) | ] \le$$
$$\le \bar{b} \, e^{\mu T} \, C [ | z_1(s) - z_2(s) | + | z_1(\sigma) - z_2(\sigma) | ]$$

This is precisely the kind of Lipschitz condition that we need for existence and uniqueness of solutions of multiple Volterra integral equations [BB] . (In the present case, the degree of multiplicity of the Volterra equation is 2.) The terms of first degree, in the Riccati – Volterra equation, clearly satisfy the Lipschitz condition for first-degree (linear) terms in a Volterra equation. ///

We present, below, an example in which Riccati-Volterra equations occur:

Example 6.1. One model of the spread of epidemics utilizes Volterra integral equations [M]. The model is an integral equation variant of the Verhulst model of population dynamics, with an additional ingredient, the so called "contact rate".
We take a model in a slightly different formulation, to conform to the notational conventions of the present paper. The model is

$$y(t) = y_0(t) + \int_0^t \; [ f_0(t,s,y(s)) + u(s) f_1(t,s,y(s)) ] \, ds$$

The state $y$ and the control $u$ are real-valued. The control $u$ represents a controllable contact rate, for instance controlled through quarantine and isolation policies.

The objective is to minimize a functional

$$J = \int_0^T \left[ \frac{\alpha}{2}(u(t))^2 + \frac{1}{2}(y(t))^2 \right] dt$$

Following the results in section… , we have the following conditions for an optimal control , the state *y*, and the co-state $\psi$ :

$$u^*(t) = -\frac{1}{\alpha}\int_t^T \psi(s) f_1(s,t,y(t))\, ds \ ;$$

$$\psi(t) = y(t) + \int_t^T \psi(s) f_{0,y}(s,t,y(t))\, ds - \frac{1}{\alpha}\int_t^T \int_t^T f_{1,y}(s,t,y(t))\, f_1(\sigma,t,y(t))\, \psi(s)\psi(\sigma)\, d\sigma\, ds \ ;$$

$$y(t) = y_0(t) + \int_0^t f_0(t,s,y(s))\, ds - \frac{1}{\alpha}\int_0^t \int_s^T f_1(t,s,y(s))\, f_1(\sigma,s,y(s))\, \psi(\sigma)\, d\sigma\, ds$$

Here, we are faced with the problem of solving simultaneously for the state *y* and the co-state $\psi$ . The equation for the co-state, in (), is of the type studied in section (), with $f_{0,y}(s,t,y(t))$ corresponding to $a(t,s)$, and $f_{1,y}(s,t,y(t))\, f_1(\sigma,t,y(t))$ corresponding to $b(t,s,\sigma)$. If we assume that the conditions of section () are met, uniformly in *y*, then there exists a constant, say $C_1$, such that $|\psi(t)| \le C_1$, uniformly in *t* and *y*. Then a Lipschitz condition for the dependence on *y* in the equation for the state will be satisfied if the function

$$F_{1,1}(t,s,y) := f_0(t,s,y) - \frac{C_1}{\alpha} f_1(t,s,y) \int_s^T f_1(\sigma,s,y)\, d\sigma$$

is Lipschitz in *y*. Under these conditions, the problem is solvable by standard iterative methods for nonlinear Volterra integral equations. ///

References.

[A]. M. de ACUTIS, "On the quadratic optimal control problem for Volterra integro-differential equations", Rend. Semin. Matem. Univ. Padova, Vol. 73, 1985, pp. 231-247.

[B]. V. L. BAKKE, "A maximum principle for optimal control problem with integral constraints", J. Optimiz. Th. Applic., Vol. 13, no. 1, 1974, pp. 32-55.

[BLT]. H. T. BANKS, B. M. LEWIS, H. T. TRAN, "Nonlinear feedback controllers and compensators: a state-dependent Riccati equation approach", Comput. Optim. Applic., Vol. 37, 2007, pp. 177-218.

[BB]. S. A. BELBAS, Y. BULKA, "Numerical solution of multiple nonlinear integral equations", Appl. Math. Comput., Vol. 217, no. 9, 2011, pp. 4791-4804.

[BS]. S. A. BELBAS, W. H. SCHMIDT, "Optimal control of Volterra integral equations with variable impulse times", Applied Mathematics and Computation, Vol. 214, 2009, pp. 353-369.

[BS1]. S. A. BELBAS, W. H. SCHMIDT, "Solvable cases of optimal control problems for integral equations", arXiv: 1606.05803, 2016.

[BLW]. S. BITTANTI, A. J. LAUB, J. C. WILLEMS, The Riccati equation, Springer-Verlag, Berlin – Heidelberg, 1991

[C]. M. A. CONNOR, "Optimal control of linear systems represented by integral equations: an iterative method", IEEE Trans. Autom. Control, Vol. 17, no. 3, 1972, pp. 404-406.

[F]. G. F. FASSHAUER, "Positive definite functionals: past, present, and future", Dolomites Research Notes on Approximation, Vol. 4, 2011, pp. 21-63.

[H]. H. HE, "On matrix-valued square integrable positive definite functions", Monatshefte f. Math., Vol. 177, no. 3, 2015, pp. 437-449.

[L]. F. L. LEWIS, Applied optimal control and estimation, Prentice-Hall, Englewood Cliffs, NJ, 1992.

[M]. J. MERCER, "Functions of positive and negative type, and their connection with the theory of integral equations", Philos. Trans. Royal Soc. London, Ser. A, Vol. 209, 1909, pp. 415-446.

[ME]. E. MESSINA, "Numerical simulation of a SIS epidemic model based on a nonlinear Volterra integral equation", Dynamical Systems, Differential Equations and Applications, AIMS Proceedings, 2015,

[P]. L. PANDOLFI, "The quadratic regulator problem and the Riccati equation for a process governed by a linear Volterra integrodifferential equation", IEEE Tans. Autom. Control, Vol. 63, no. 5, 2018, pp. 1517-1522,

[PY]. A. J. PRITCHARD, Y. YOU, “Causal feedback optimal control for Volterra integral equations”, SIAM J. Control Optimiz., Vol. 34, 1996, pp. 1874-1896.

[S1]. W. H. SCHMIDT, “Notwendige Optimalitätsbedingungen für Prozesse mit zeitvariablen Integralgleichungen in Banach-Räumen”, Z. angew. Mathematik u. Mechanik, Vol. 60, 1980, pp. 595-608.

[S2]. W. H. SCHMIDT, “Durch Integralgleichungen beschriebene Optimale Prozesse mit Nebenbedingungen in Banachräumen – notwendige Optimalitätsbedingungen”, Z. angew. Mathematik u. Mechanik, Vol. 62, 1982, pp. 65-75.

[W]. W. von WALDENFELS, “ Positiv definite Funktionen auf einem unendlichdimensionalen Vektorraum”, Studia Mathematica, Vol. 30, 1968, pp. 153-162.